\documentclass[a4paper,12pt,reqno]{amsart}

\usepackage{amsmath}
\usepackage{amsthm}
\usepackage{amssymb}
\usepackage{graphicx}
\usepackage{appendix}
\usepackage[pdftex]{hyperref}
\hypersetup{breaklinks=true,colorlinks=true,citecolor=green,urlcolor=red,linkcolor=blue,pdfpagemode=UseNone}

\theoremstyle{plain}

\newtheorem{theorem}{Theorem}[section]

\newcommand{\eps}{\varepsilon}

\renewcommand{\leq}{\leqslant}

\def\COMMENT#1{}
\let\COMMENT=\footnote

\def\E{\mathbb{E}}
\DeclareMathOperator{\Var}{Var}
\DeclareMathOperator{\Aut}{Aut}
\DeclareMathOperator{\GL}{GL}
\DeclareMathOperator{\tr}{tr}
\DeclareMathOperator{\cl}{cl}
\DeclareMathOperator{\Cl}{Cl}

\def\Markstrom{Markstr\"om}
\def\Umea{Ume\aa}
\def\Malostranske{Malostransk\'{e}}
\def\Namesti{N\'am\v{e}st\'{i}}
\def\Frobenius{Fr\"obenius}

\begin{document}

\title{The thresholds for diameter 2 in random Cayley graphs}
\author{Demetres Christofides}
\author{Klas \Markstrom}
\thanks{During the work of this project DC was funded from the European Research Council under the European Union's Seventh Framework Programme (FP7/2007-2013)/ERC grant agreement no. 259385.} 
\date{\today}
\begin{abstract}
Given a group $G$, the model $\mathcal{G}(G,p)$ denotes the probability space of all Cayley graphs of $G$ where each element of the generating set is chosen independently at random with probability $p$.

In this article we show that for any $\varepsilon > 0$ and any family of groups $G_k$ of order $n_k$ for which $n_k \to \infty$, a graph $\Gamma_k \in \mathcal{G}(G_k,p)$ with high probability has diameter at most 2 if $p \geqslant \sqrt{(2 + \eps) \frac{\log{n_k}}{n_k}}$ and with high probability has diameter greater than 2 if  $p \leqslant \sqrt{\left(\frac{1}{4} + \eps\right)\frac{\log{n_k}}{n_k}}$. 

We also provide examples of families of graphs which show that both of these results are best possible.

Of particular interest is that for some families of groups, the corresponding random Cayley graphs achieve diameter~2 significantly faster than the Erd\H{o}s-Renyi random graphs.
\end{abstract}
\maketitle

\section{Introduction}

Let us begin by recalling that given a group $G$ and a subset $S$ of $G$, the Cayley graph $\Gamma = \Gamma(G;S)$ of $G$ with respect to $S$ has the elements of $G$ as its vertex set and has an edge between $g$ and $h$ if and only if $hg^{-1} \in S$ or $gh^{-1} \in S$. We ignore any loops or multiple edges. In particular, whether $1 \in S$ or not is immaterial. Observe for example that $\Gamma$ is connected if and only if the set $S$ generates the group $G$. Throughout the paper, we will often refer to the set $S$ as the generating set of the graph $\Gamma$ irrespectively of whether it is a generating set for the group $G$ or not. 

The model $\mathcal{G}(G,p)$ is the probability space of all graphs $\Gamma(G;S)$ in which every element of $G$ is assigned to the set $S$ independently at random with probability $p$. This model has many similarities with the model $\mathcal{G}(n,p)$, which is the probability space of all graphs with vertex set $\{1,2,\ldots,n\}$ in which every edge appears independently with probability $p$. We refer the reader to~\cite{CM} for some of these similarities. There are however many differences between these two models. An obvious difference is that every graph in $\mathcal{G}(G,p)$ is regular while with high probability this is not the case in the model $\mathcal{G}(n,p)$ (unless in the trivial cases in which $p$ is either so large or so small that with high probability forces $G$ to be complete or empty respectively). This difference also motivates the comparison of the model $\mathcal{G}(G,p)$ with the model $\mathcal{G}_{n,r}$, the probability space of all $r$-regular graphs on $\{1,2,\ldots,n\}$ taken with the uniform measure. These two models still have significant differences. For example, every graph $\Gamma \in \mathcal{G}(G,p)$ is not only regular, but in fact it has a high degree of symmetry. More specifically, every element $g$ of $G$ defines an automorphism of $\Gamma$ by right multiplication and so $G$ is a subgroup of $\Aut(\Gamma)$. On the other hand, it is known that for every $3 \leqslant r \leqslant n-4$, graphs in $\mathcal{G}_{n,r}$ have with high probability a trivial automorphism group~\cite{KSV}. Given the success of random graphs in settling many graph theory questions and the fact that random Cayley graphs have some important properties not shared by other random graph models, we see that their study is highly desirable. 

In this paper we study the diameter of random Cayley graphs. In~\cite{CM} it was proved that if $c > 1$ is a fixed constant, then every $\Gamma \in \mathcal{G}(G,c\log_2{|G|})$ has with high probability logarithmic diameter. On the other hand, if $c \leqslant 1$ then there are groups $G$ for which $\Gamma$ will with high probability be disconnected. Here we will be concerned with a different range of $p$. More specifically, we will be concerned with the range of $p$ for which the diameter of the random Cayley graph becomes larger than 2. So let us begin by reviewing what is know in the $\mathcal{G}(n,p)$ and $\mathcal{G}_{n,r}$ models. It is well-known~\cite{MM} that for any $\eps > 0$, a graph from $\mathcal{G}(n,p)$ with high probability has diameter at most 2 if $p \geqslant \sqrt{\frac{(2 + \eps) \log{n}}{n}}$ and diameter greater than 2 if $p \leqslant \sqrt{\frac{(2 - \eps) \log{n}}{n}}$. For the case of random regular graphs, it seems that the threshold for diameter~2 is still not known. If the sandwich conjecture of Kim and Vu~\cite{KV} is true then the expected result should hold. Namely, that if $r \leqslant \sqrt{(2 - \eps)n\log{n}}$ then almost every $r$-regular graph on $n$ vertices has diameter greater than 2, while if $r \geqslant \sqrt{(2 + \eps)n\log{n}}$, then almost every $r$-regular graph on $n$ vertices has diameter equal to 2. 

Note that when speaking about random Cayley graphs, the group $G$ and thus its order is fixed. So strictly speaking it does not really make sense to speak of events occurring with high probability. Instead, we have to consider a family of groups $G_k$ for which their orders tends to infinity. Motivated by the above results it is natural to conjecture that for any $\eps > 0$, if $(G_k)$ is a family of groups having order $n_k$ with $n_k \to \infty$ then there is a constant $c$ such that if $\Gamma \in \mathcal{G}(G_k,p)$ then with high probability $G$ has diameter greater than 2 if $p \leqslant  \sqrt{(c - \eps)\frac{\log{n}}{n}}$ and diameter at most 2 if $p \geqslant  \sqrt{(c + \eps)\frac{\log{n}}{n}}$. It turns out that this conjecture is wrong for the following reason: We will see that the conjecture is true for several natural families of groups. However, different families can give rise to different constants $c$. In particular, if we take two families $(G_k)$ and $(H_k)$ for which the conjecture is true but with different values of $c$, then if we interlace these families to a new one we get a family for which this naive conjecture is false. This suggests that the right conjecture to consider is the following: For any $\eps > 0$, if $(G_k)$ is a family of groups having order $n_k$ with $n_k \to \infty$ then there are constant $c_1,c_2$ such that if $\Gamma \in \mathcal{G}(G_k,p)$ then with high probability $G$ has diameter greater than 2 if $p \leqslant  \sqrt{(c_1 - \eps)\frac{\log{n}}{n}}$ and diameter at most 2 if $p \geqslant  \sqrt{(c_2 + \eps)\frac{\log{n}}{n}}$. Our aim in this paper is to prove this conjecture. Moreover we will show that the values we obtain for $c_1$ and $c_2$ are best possible. More specifically we will prove the following results:

\begin{theorem}\label{T:D2}
Let $\eps > 0$, let $G$ be a group on $n$ elements, and let $\Gamma \in \mathcal{G}(G,p)$, where $p \geqslant \sqrt{(2 + \eps)\frac{\log{n}}{n}}$. Then the diameter of $\Gamma$ is with high probability at most $2$.
\end{theorem}

As we have already mentioned, it does not make sense to speak about asymptotic results for a fixed group. However throughout the paper we choose to abuse notation. The above result should be interpreted as a result for families of groups instead. That this result is best possible is shown by the following:

\begin{theorem}\label{E:D>=3}
Let $\eps > 0$ and let $\Gamma \in \mathcal{G}(\mathbb{Z}_2^n,p)$, where $p \leqslant \sqrt{(2 - \eps)\frac{\log{N}}{N}}$ and $N = 2^n = |\mathbb{Z}_2^n|$. Then with high probability, the diameter of $\Gamma$ is greater than 2.
\end{theorem}

In the other direction we have the following results:

\begin{theorem}\label{T:D>=3}
Let $\eps > 0$, let $G$ be any group of order $n$ and let $\Gamma \in \mathcal{G}(G,p)$, where $p \leqslant \sqrt{\left(\frac{1}{4} - \eps \right)\frac{\log{n}}{n}}$. Then with high probability, the diameter of $\Gamma$ is greater than 2.
\end{theorem}

For abelian groups we can actually prove a better bound:

\begin{theorem}\label{T:D>=3-Abelian}
Let $\eps > 0$, let $G$ be any abelian group of order $n$ and let $\Gamma \in \mathcal{G}(G,p)$, where $p \leqslant \sqrt{\left(\frac{1}{2} - \eps \right)\frac{\log{n}}{n}}$. Then with high probability, the diameter of $\Gamma$ is greater than 2.
\end{theorem}

Theorems~\ref{T:D>=3} and~\ref{T:D>=3-Abelian} are also best possible. We begin by stating the corresponding result for abelian groups.

\begin{theorem}\label{E:D>=3-Abelian}
Let $\eps > 0$, let $G$ be the cyclic group of order $n$ and let $\Gamma \in \mathcal{G}(G,p)$, where $p \geqslant \sqrt{\left(\frac{1}{2} + \eps \right)\frac{\log{n}}{n}}$. Then with high probability, the diameter of $\Gamma$ is at most 2.
\end{theorem}

Before stating the corresponding result for general groups we recall that for an element $x$ of a group $G$ its conjugacy class is the set $\Cl(x) = \{y^{-1}xy:y \in G\}$ and its size is denoted by $\cl(x)$. We also recall that an element $x$ of $G$ with $x^2=1$ is called an involution.

\begin{theorem}\label{E:D=2}
Let $0 < \eps < 1/4$ and let $G$ be a group of order $n$ such that
\begin{itemize}
\item[(a)] $G$ contains at most $O(n^{(1 + \eps)/2})$ involutions.
\item[(b)] $G$ contains at most $O(n^{(1 + \eps)/2})$ elements $x$ with $\cl(x) \leqslant 1/\eps$.
\item[(c)] $G$ contains at most $O(n^{(1 + \eps)/4})$ involutions $x$ with $\cl(x) \leqslant 1/\eps$.
\end{itemize}
Let also $\Gamma \in \mathcal{G}(G,p)$, where $p \geqslant \sqrt{\left(\frac{1}{4} + \eps \right) \frac{\log{n}}{n}}$. Then with high probability, the diameter of $\Gamma$ is at most 2.
\end{theorem}

Observe that for any such group, if we choose $p = \sqrt{\left(\frac{1}{4} + \eps \right) \frac{\log{n}}{n}}$, then the corresponding Cayley graph has with high probability at most $\sqrt{\left(\frac{1}{4} + 2\eps\right)n^3\log{n}}$ edges and yet has diameter~2. This is in contrast with the Erd\H{o}s-Renyi random graph case in which we need about $\sqrt{2}$ times as many edges in order for the diameter to become equal to~2.

Of course, in order to show that Theorem~\ref{T:D>=3} is best possible, we need to actually exhibit a family of groups satisfying properties (a)-(c) of Theorem~\ref{E:D=2}. One can easily check that the symmetric groups satisfy these properties. Indeed the smallest non-trivial conjugacy class of $S_n$ is the conjugacy class of $(1 \, 2)$ which contains $\binom{n}{2}$ elements. So (b) and (c) trivially hold. To see that (a) holds observe that if $a_n$ is the number of involutions of $S_n$ then it satisfies the recurrence relation $a_n = a_{n-1} + (n-1)a_{n-2}$. Even though one can give finer asymptotics for the number of involutions (see e.g.~\cite[Example~5.17]{Wilf}) it is easy to prove by induction using the recurrence relation that $a_n \leqslant 2^n (n!)^{1/2} = (n!)^{1/2 + n\log{2}/\log{(n!)}}$ which implies (a). 

In~\cite{CM} we introduced a model of random graphs based on Latin squares. Recall that a Latin square of order $n$ is an $n \times n$ matrix $L$ with entries from a set of $n$ elements, such that in each row and in each column, every element appears exactly once. Given a Latin square $L$ with entries in a set $A$ of size $n$, and a subset $S$ of $A$, we define the Latin square graph $\Gamma(L;S)$ on vertex set $[n]$, by joining $i$ to $j$ if and only if either $L_{ij} \in S$ or $L_{ji} \in S$. Observe that with this definition every Cayley graph can be seen as a Latin square graph by taking the rows and columns of $L$ to be indexed by the elements of $G$ and defining $L_{xy} = xy^{-1}$. Then for any subset $S$ of $G$ we have that $\Gamma(G;S)$ is exactly the same as $\Gamma(L;S)$. The model $\mathcal{G}(L,p)$ of random Latin square graphs is defined with exact analogy to the model $\mathcal{G}(G,p)$. Namely it is the probability space of all graphs $\Gamma(L; S)$ in which every element from the set of entries of $L$ is assigned to the set $S$ independently at random with probability $p$.

It is natural therefore to ask how the above results generalise to random Latin square graphs. We will prove the following result.

\begin{theorem}\label{Latin}
Let $\eps > 0$, let $L$ be any latin square of order $n$ and let $\Gamma \in \mathcal{G}(L,p)$. If $p \geqslant \sqrt{\left(26 + \eps \right)\frac{\log{n}}{n}}$ then with high probability the diameter of $\Gamma$ is at most~2, while if $p \leqslant \sqrt{\left(\frac{1}{4} - \eps \right)\frac{\log{n}}{n}}$ then with high probability, the diameter of $\Gamma$ is greater than 2.

\end{theorem}



We stated all of the above results with $\eps$ being a fixed positive constant. In fact it will be clear from the proofs that $\eps$ can be allowed to vary with $n$. In fact all the results will still hold as long as $\eps = \omega(1/\log{n})$.

We now give an overview of the structure of the paper. In Section~2 we collect the probabilistic tools and in Section~3 we collect the representation theoretic results that we will use. In Section~4 we prove Theorems~\ref{E:D>=3},~\ref{T:D>=3} and~\ref{T:D>=3-Abelian}. Section~5 is devoted to the proofs of Theorems~\ref{T:D2} and~\ref{E:D>=3-Abelian}. In Section~6 we prove Theorem~\ref{E:D=2}. Section~6 is the only place where the representation theoretic results will be used. Sections~4 and~5 can be read independently while in Section~6 we will use the notation and ideas introduced in Section~5. Finally in Section~7 we use ideas from Sections~4 and~5 to prove Theorem~\ref{Latin}.

\section{Probabilistic Tools}

In the proofs of Theorems~\ref{T:D2}-\ref{E:D=2} we will make use of the following probabilistic tools. We refer the reader to~\cite{AS} for their proofs.

\begin{theorem}[Chebyshev's Inequality]
Let $X$ be a random variable which takes values on the non-negative integers and suppose that it has finite variance and expectation. Then
\[
\Pr(X = 0) \leqslant \frac{\Var(X)}{(\E X)^2}.
\]
\end{theorem}

\begin{theorem}[Kleitman's Inequality]
Let $\Omega$ be a finite set and let $\{F_i\}_{i \in I}$ be subsets of $\Omega$, where $I$ is a finite index set. Let $R$ be a random subset of $\Omega$ and for each $i \in I$ let $E_i$ be the event that $F_i \subseteq R$. Then 
\[
\Pr\left(\bigcap_{i\in I} \overline{E_i} \right) \geqslant \prod_{i \in I} \Pr(\overline{E_i}).
\]
\end{theorem}

\begin{theorem}[Janson's Inequality]
Let $\Omega$ be a finite set and let $\{F_i\}_{i \in I}$ be subsets of $\Omega$, where $I$ is a finite index set. Let $R$ be a random subset of $\Omega$ and for each $i \in I$ let $E_i$ be the event that $F_i \subseteq R$. Suppose also that $\Pr(E_i) \leqslant \eps$ for each $i \in I$. Then
\[
\Pr\left(\bigcap_{i\in I} \overline{E_i} \right) \leqslant \exp\left(-\sum_{i \in I} \Pr(E_i) + \sum_{i \in I}\sum_{\{j \neq i : F_i \cap F_j \neq \emptyset\}} \Pr(E_i \cap E_j) \right).
\]
\end{theorem}

\section{Representation Theoretic Tools}

In this section we recall several results from representation theory that we will need for the proof of Theorem~\ref{E:D=2}. The proofs can be found in many books on representation theory, for example in~\cite{JL}. Throughout this section we will assume that all groups are finite.

Given two elements $x,y$ of a group $G$ we say that $x$ is conjugate to $y$ if there is an element $z$ of $G$ such that $y = z^{-1}xz$. This is easily seen to be an equivalence relation and the equivalence classes of this relation are called the conjugacy classes. We will denote the number of equivalence classes of $G$ by $\cl(G)$.

A representation $\rho$ of a group $G$ is a homomorphism $\rho : G \to \GL(V)$, where $V$ is a finite dimensional vector space over $\mathbb{C}$. The dimension $d_{\rho}$ of $\rho$ is simply the dimension of $V$. We say that a subspace $W$ of $V$ is invariant if it is fixed by $\rho$. (I.e.~$\rho(g)(W) \subseteq W$ for every $g \in G$.) It is then easily checked that the restriction $\rho_W : G \to \GL(W)$ is a representation. We say that $\rho$ is irreducible if there is no non-trivial invariant subspace.

Given a representation $\rho$ of a group $G$, its character is the function $\chi_{\rho}:G \to \mathbb{C}$ defined by $\chi_{\rho}(g) = \tr(\rho(g))$, where, $\tr(A)$ denotes the trace of the linear transformation $A$. We say that the character is irreducible if the corresponding representation is irreducible. Every group has only finitely many irreducible characters. In fact the following result holds.

\begin{theorem}\label{dimensions}
Let $G$ be a group of order $n$ and let $R$ be the set of all irreducible representations of $G$. Then $|R| = \cl(G)$ and moreover
\[ \sum_{\rho \in R} d_{\rho}^2 = n.\]
\end{theorem}



Given an irreducible character $\chi$ of a group $G$, its \Frobenius-Schur indicator is defined by
\[
\iota(\chi) = \frac{1}{|G|}\sum_{g \in G} \chi(g^2).
\]
We will use the following two properties of the \Frobenius-Schur indicator

\begin{theorem}\label{Frobenius-Schur} $ $
\begin{itemize}
\item[(a)] If $\chi$ is an irreducible character of a group $G$, then its \Frobenius-Schur indicator $\iota(\chi)$ takes values in $\{-1,0,1\}$.
\item[(b)] If $g$ is an element of a group $G$ then
\[ \sum_{\chi} \iota(\chi) \chi(g) = |\{h \in G: h^2 = g\}|,\]
where the sum is over all irreducible representations of $G$.
\end{itemize}
\end{theorem}

The following is the main result of this section that we will need in our proof of Theorem~\ref{E:D=2}.

\begin{theorem}\label{main-rep}
Let $G$ be a group of order $n$ and let $x$ be an element of $G$. Then there are at most $\sqrt{n\cl(G)}$ elements $y$ of $G$ such that $y^2 = x$. In particular, all but at most $7\sqrt{n\cl(G)}$ elements $y$ of $G$ satisfy $y^2,(y^{-1}x)^2 \notin \{1,x,x^{-1},x^2\}$, where 1 denotes the identity element of $G$. 
\end{theorem}

\begin{proof}
It is enough to prove only the first claim as the last one follows immediately from it after observing that $y^2 = x$ if and only if $(y^{-1}x)^2 = x$. By part (b) of Theorem~\ref{Frobenius-Schur}, the number of $y$ such that $y^2 = x$ is equal to $\sum_{\chi} \iota(\chi) \chi(x)$ and by part (a) of Theorem~\ref{Frobenius-Schur} this is at most $\sum_{\chi} |\chi(x)|$. Given an irreducible character $\chi$ of $G$, let $\rho = \rho_{\chi}$ be the corresponding linear transformation. Since $\rho(x)^n$ is the identity matrix, it follows that all eigenvalues of $\rho(x)$ are $n$-th roots of unity and so $|\chi(x)| = |\tr(\rho(x))| \leqslant d_{\rho}$. In particular, by Theorem~\ref{dimensions} and the Cauchy-Schwarz inequality we get that the number of elements which square to $x$ is at most
\[
\sum_{\rho} d_{\rho} \leqslant \left(\sum_{\rho} d_{\rho}^2 \right)^{1/2} \left(\sum_{\rho} 1 \right)^{1/2} = \sqrt{n\cl(G)},
\]
where all sums above are over all irreducible representations of $G$.
\end{proof}

We will also need the following result.

\begin{theorem}\label{commuting_pairs}
Let $G$ be a group of order $n$. Then for every element $x$ of $G$
\begin{itemize}
\item[(a)] there are (at most) $n/\cl(x)$ elements $y$ such that $y^{-1}xy=x$ and
\item[(b)] there are at most $n/\cl(x)$ elements $y$ such that $y^{-1}xy=x^{-1}$.
\end{itemize}
\end{theorem}
 
\begin{proof}
Indeed part (a) is an immediate consequence of the orbit-stabiliser theorem. For part (b), we just observe that for each $x$ there will be either $n/\cl(x)$ or 0 elements $y$ such that $y^{-1}xy = x^{-1}$ depending on whether $x$ and $x^{-1}$ belong to the same conjugacy class or not.   
\end{proof} 
 
\section{Proofs of Theorems~\ref{E:D>=3},\ref{T:D>=3} and~\ref{T:D>=3-Abelian}}

We begin with the proof of Theorem~\ref{E:D>=3} which is a simple application of the second moment method.

\begin{proof}[Proof of Theorem~\ref{E:D>=3}]
Let us write $0$ for the identity element of $\mathbb{Z}_2^n$ and let $X$ be the number of vertices $x$ of $\Gamma$ which are at distance greater than two from $0$. It suffices to prove that with high probability $X \neq 0$. Observe that the distance between 0 and $x$ is greater than 2 if and only if $x$ does not belong to the generating set and moreover there is no pair of the form $\{y,y+x\}$ with $y \neq 0,x$ for which both of its elements belong to the generating set. Since there are exactly $(N-2)/2$ such pairs which are pairwise disjoint, we have that
\[ \E X = (N - 1)(1 - p)(1 - p^2)^{(N-2)/2}. \]
We also claim that 
\[
\E X(X-1) = (N-1)(N-2)(1-p)^2 \left( (1 - p)^4 + 4p(1-p)^3 + 2p^2(1-p^2) \right)^{(N-4)/4}. 
\]
To see this, given distinct vertices $x,y$ of $\Gamma \setminus \{0\}$, partition the elements of $G \setminus \{0,x,y,x+y\}$ into quadruples of the form $\{z,z+x,z+y,z+x+y\}$ and observe that both $x$ and $y$ are at distance greater than 2 from 0 if and only if $x,y$ do not belong to the generating set and moreover for each quadruple of the form $\{z,z+x,z+y,z+x+y\}$, none of its subsets of the form $\{z,z+x\},\{z,z+y\},\{z+x,z+x+y\},\{z+y,z+x+y\}$ is contained in the generating set.

Since $p \leq \sqrt{(2-\eps) \frac{\log{N}}{N}}$, we have that 
\[
\E X = (1 + o(1))N\exp\left\{-p^2N/2 + o(1) \right\} = \Omega(N^{\eps/2})
\]
and
\[
\E X(X-1) = (1 + o(1))N^2\exp\left\{-p^2N + o(1) \right\}.
\]
Thus by Chebyshev's inequality we have
\[
\Pr(X = 0) \leqslant \frac{\Var(X)}{(\E X)^2} = \frac{\E X(X-1)}{(\E X)^2} - 1 + \frac{1}{\E X} = o(1),
\]
as required.  
\end{proof}

We would like to follow the same approach as above in order to prove Theorem~\ref{T:D>=3}. In fact, it is relatively straightforward to show that if $\Gamma \in \mathcal{G}(G,p)$, where $p \leqslant \sqrt{\left(\frac{1}{4} - \eps\right)\frac{\log{n}}{n}}$ and $X$ is the number of vertices of $\Gamma$ which are at distance greater than two from the identity of $G$ then $\E X$ tends to infinity. It seems though that without knowing anything about the group structure it is difficult to work out good approximations of $\mathbb{E}X$ and $\Var(X)$. We can only find upper and lower estimates for them and unfortunately it seems difficult to show directly that $\Var(X)/\E X^2 \to 0$. For this reason we will concentrate not on $X$ but on a related random variable $Y$ which we can have better control over. 

\begin{proof}[Proof of Theorem~\ref{T:D>=3}]

Let us write $1$ for the identity element of $G$. For each vertex $x \neq 1$ of $\Gamma$, we define a graph $\Gamma_x$ whose vertices are the elements of $G$ and where there is an edge between $g$ and $h$ if and only if there is a $y$ such that the appearance of $g$ in the generating set guarantees an edge between $1$ and $y$ and the appearance of $h$ guarantees an edge between $y$ and $x$ or vice versa. (I.e.~$h$ guarantees an edge between $1$ and $y$ and $g$ between $y$ and $x$.) Note that we do allow $\Gamma_x$ to have loops. An easy computation shows that $h$ is adjacent to $g$ in $\Gamma_x$ if and only if $h \in \{xg,xg^{-1},x^{-1}g,x^{-1}g^{-1},gx,gx^{-1},g^{-1}x,g^{-1}x^{-1}\}$. In particular, $\Gamma_x$ has maximum degree at most $8$ and so it contains at most $4n$ edges. It is also immediate that $g$ is adjacent to $gx$ in $\Gamma_x$ and so $\Gamma_x$ has minimum degree at least $1$ and therefore it contains at least $n/2$ edges. Note also that $g$ and $h$ are adjacent in $\Gamma_{x}$ if only if $x \in \{hg,hg^{-1},h^{-1}g,h^{-1}g^{-1},gh,gh^{-1},g^{-1}h,g^{-1}h^{-1}\}$. Thus for every two distinct elements $g,h$ of $G$ there are at most eight distinct $\Gamma_x$'s in which $g$ and $h$ are adjacent.

We now claim that $G$ contains a set $A$ of size at least $n^{1-\eps}$ which does not contain the identity of $G$ and such that for any distinct elements $x,y$ of $A$ we have that $\Gamma_x$ and $\Gamma_y$ have at most $30n^{1-\eps}$ common edges. Indeed let us create a new graph $H$ whose vertices are the elements of $G$ excluding the identity and where $x$ is adjacent to $y$ if and only if $\Gamma_x$ and $\Gamma_y$ have more than $30n^{1-\eps}$ common edges. Since $\Gamma_x$ has at most $4n$ edges and every edge of $\Gamma_x$ belongs to at most 7 other $G_y$'s we deduce that the maximum degree in $H$ is at most $14n^{\eps}/15$. So the independence number of $H$ is at least $15n^{1-\eps}/14$ and so the claim follows.

Now let $Y$ be the number of vertices $x$ of $A$ which are at distance greater than two from the identity. It suffices to show that with high probability $Y \neq 0$. We begin by calculating the expectation of $Y$. For any element $x$ of $G$, let us write $B_x$ for the event that $x$ has distance greater than two from the identity. Then we have that
\[
\Pr(B_x) \geqslant (1-p)^2(1-p^2)^{4n}.
\]
This follows directly from Kleitman's Inequality since $B_x = \overline{E_x} \cap \overline{E_{x^{-1}}} \bigcap_{e \in E(\Gamma_x)} \overline{E_e}$, where $E_x,E_{x^{-1}}$ denote the events that $x,x^{-1}$ appear to the generating set respectively and for an edge $e$ of $\Gamma_x$, $E_e$ denotes the event that both of its incident vertices appear in the generating set.  

It follows that
\[
\E Y \geqslant n^{1-\eps}(1-p)^2(1-p^2)^{4n} = n^{1-\eps}\exp\left\{-4p^2n + o(1)\right\} = \Omega(n^{3\eps}).
\]
Our next aim is to show that $\E Y(Y-1)$ is asymptotically equal to $(\E Y)^2$. We have that
\[
\E Y(Y-1) = \sum_{x\in Y}\sum_{y \in Y \setminus{x}} P(B_x \cap B_y). 
\]
Recall that $\Gamma_x$ and $\Gamma_y$ have at most $30n^{1-\eps}$ common edges. Let us write $e_1,\ldots,e_r$ for the common edges of $\Gamma_x$ and $\Gamma_y$ and let $e_{r+1},\ldots,e_s$ be all other edges of $\Gamma_y$. For each $1 \leqslant i \leqslant s$, let $C_i$ be the event that not both vertices incident to $e_i$ appear in the generating set and let $D_i = \cap_{j=1}^i C_j$. Since $D_s = B_y$ then
\[
\Pr(B_x \cap B_y) = \Pr(B_x) \Pr(B_y) \frac{\Pr(B_x|B_y)}{\Pr(B_x)} = \Pr(B_x) \Pr(B_y) \prod_{i=1}^s \frac{\Pr(B_x|D_i)}{\Pr(B_x|D_{i-1})}
\]
where by convention $\Pr(B_x|D_0) := \Pr(B_x)$. 

By the law of total probability, for each $1 \leqslant i \leqslant s$ we have
\begin{align*}
\Pr(B_x|D_{i-1}) &= \frac{\Pr(B_x \cap D_{i-1})}{\Pr(D_{i-1})} = \frac{\Pr(B_x \cap D_{i-1} \cap C_i) + \Pr(B_x \cap D_{i-1} \cap \overline{C_i})}{\Pr(D_{i-1})} \\
&= \frac{\Pr(B_x | D_{i-1} \cap C_i)\Pr(D_{i-1} \cap C_i)}{\Pr(D_{i-1})} + \frac{\Pr(B_x \cap D_{i-1} \cap \overline{C_i})}{\Pr(D_{i-1})} \\
&= \Pr(B_x|D_i) \Pr(C_i|D_{i-1}) + \frac{\Pr(B_x \cap D_{i-1} \cap \overline{C_i})}{\Pr(D_{i-1})}.
\end{align*}

To bound the first term, we use Kleitman's Lemma to the events $\overline{C_i}$ and $\overline{D_{i-1}}$ to get 
\[
\Pr(B_x|D_i) \Pr(C_i|D_{i-1}) \geqslant \Pr(B_x|D_i)\Pr(C_i) = (1-p^2)\Pr(B_x|D_i).
\]
To bound the second term, observe first that if $1 \leqslant i \leqslant r$ then $B_x \subseteq C_i$ and so $\Pr(B_x \cap D_{i-1} \cap \overline{C_i}) = 0$. If $r+1 \leqslant i \leqslant s$, let $X_i$ denote the event that no vertex incident to a vertex of $e_i$ in $\Gamma_x$ appears in the generating set and $Y_i$ the event that for each $j < i$ if $e_j$ meets $e_i$ then the vertex incident to $e_j$ but not $e_i$ is not in the generating set. Then, $B_x \cap D_{i-1} \cap \overline{C_i} = B_x \cap D_{i-1} \cap \overline{C_i} \cap X_i \cap Y_i$. Observe now that whether the two vertices incident to the edge $e_i$ appear in the generating set or not does not affect the outcome of the event $B_x \cap D_{i-1} \cap X_i \cap Y_i$. Applying now Kleitman's lemma to the events $\overline{B_x} \cup \overline{D_{i-1}},\overline{X_i}$ and $\overline{Y_i}$ we get
\[
\Pr(B_x \cap D_{i-1} \cap \overline{C_i}) \geqslant \Pr(\overline{C_i})\Pr(X_i)\Pr(Y_i)\Pr(B_x \cap D_{i-1}) \geqslant p^2(1-p)^{28}\Pr(B_x|D_{i-i})\Pr(D_{i-1}).
\]
To see the last inequality, recall that $\Gamma_x$ and $\Gamma_y$ have maximum degree at most $8$ and so each of the events $X_i$ and $Y_i$ says that at most 14 elements of $G$ (at most seven for each vertex incident to $e_i$) do not appear in the generating set.

It follows that if $1 \leqslant i \leqslant r$ then $\Pr(B_x|D_i) \leqslant \Pr(B_x|D_{i-1})/(1-p^2)$ while if $r+1 \leqslant i \leqslant s$ then
\[
\frac{\Pr(B_x|D_i)}{\Pr(B_x|D_{i-1})} \leqslant \frac{1 - p^2(1-p)^{28}}{1 - p^2} \leqslant 1 + 29p^3
\]
provided that $n$ is large enough. So putting everything together we get
\[
\Pr(B_x \cap B_y) \leqslant \Pr(B_x) \Pr(B_y) (1-p^2)^{-30n^{1-\eps}}(1 + 29p^3)^{4n} \leqslant (1 + o(1))\Pr(B_x) \Pr(B_y).
\]
Thus $\E Y(Y-1) \leqslant (1 + o(1))(\E Y)^2$ and thus by Chebyshev's inequality,
\[
\Pr(Y = 0) \leqslant \frac{\Var(Y)}{(\E Y)^2} = \frac{\E Y(Y-1)}{(\E Y)^2} - 1 + \frac{1}{\E Y} = o(1).
\]
This completes the proof.
\end{proof}

For specific families of groups we may be able to do better. For example if the group $G$ is abelian, then for each non-identity element $x$ of $G$, the graph $\Gamma_x$ defined in the proof of Theorem~\ref{T:D>=3} has maximum degree four as $g$ can only be adjacent to $xg,xg^{-1},x^{-1}g$ and $x^{-1}g^{-1}$ in $\Gamma_x$. It is now straightforward to modify this proof to deduce Theorem~\ref{T:D>=3-Abelian}.

\section{Proofs of Theorems~\ref{T:D2} and~\ref{E:D>=3-Abelian}}

To show the proof idea we begin by proving a similar result but with a worse constant.

\begin{theorem}\label{T:D2-worse constant}
Let $\eps > 0$, let $G$ be a group on $n$ elements, and let $\Gamma \in \mathcal{G}(G,p)$, where $p \geqslant \sqrt{(7+\eps)\frac{\log{n}}{n}}$. Then the diameter of $\Gamma$ is with high probability at most $2$. 
\end{theorem}

\begin{proof}
Let $1$ denote the identity element of $\Gamma$. Since $\Gamma$ is vertex transitive it is enough to show that with high probability, for every $x \in G$, $d_{\Gamma}(1,x) \leqslant 2$. We claim that for every $x \in G$ we have 
\[\Pr(d(1,x) > 2) \leqslant (1 - p^2)^{(n-2)/7}.\]
The result will then follow since the probability that there is an $x \in G$ with $d_{\Gamma}(1,x) > 2$ is at most 
\[
n(1 - p^2)^{(n-2)/7} \leqslant n\exp\left\{-\frac{p^2(n-2)}{7}\right\} =  o(1).
\]
To prove the claim, let us fix an $x \neq 1$ and for each $y \neq 1,x$ let us denote by $A_y := A_y(1,x)$ the event that the edges between $1$ and $y$ and between $x$ and $y$ both appear in $\Gamma$. So, if $d(1,x) > 2$ then none of the events $A_y$ occurs and so 
\[ \Pr(d(1,x) > 2) \leqslant \Pr \left( \bigcap_{y \neq 1,x} \overline{A_y} \right).\]
Note that for each $y$ we have that $\Pr(A_y) \geqslant p^2$ since for example the occurrence of $y$ and $yx^{-1}$ in the generating set guarantees that the event $A_y$ occurs. In particular, we have that $\Pr(\overline{A_y}) \leqslant 1 - p^2$. The events $\{A_y : y \in G - \{1,x\} \}$ are not necessarily independent. We claim however that there is a subset $I$ of $G \setminus \{1,x\}$ with $|I| \geqslant (n-2)/7$ such that the events $\{A_y:y \in I\}$ are independent. Our earlier claim will then follow as then we will have
\[ 
\Pr(d(1,x) > 2) \leqslant \Pr \left( \bigcap_{y \neq 1,x} \overline{A_y} \right) \leqslant \Pr \left( \bigcap_{y \in I} \overline{A_y} \right) \leqslant (1-p^2)^{(n-2)/7}. 
\]  
It remains therefore to find such a set $I$. Observe that the event $A_y$ depends only on the occurrence of the elements $y,y^{-1},yx^{-1}$ and $xy^{-1}$ in the generating set. Consider the natural dependency graph $H$ on vertex set $G \setminus \{1,x\}$ in which $y$ is adjacent to $z$ if and only if
\[
\{y,y^{-1},yx^{-1},xy^{-1}\} \cap \{z,z^{-1},zx^{-1},xz^{-1}\} \neq \emptyset.
\]
Then, if $I$ is any independent set of $H$, the events $\{A_y:y \in I\}$ are mutually independent. So to complete the proof of the claim and thus the proof of the theorem it is enough to show that $H$ has an independent set of size at least $(n-2)/7$. But $H$ has maximum degree six as the only possible neighbours of $y$ in $H$ are $y^{-1},yx^{-1},xy^{-1},xy^{-1}x,yx$ and $y^{-1}x$. These are indeed the only choices for $z$ for which 
\[
\{y,y^{-1},yx^{-1},xy^{-1}\} \cap \{z,z^{-1},zx^{-1},xz^{-1}\} \neq \emptyset.
\] 
Since $H$ has maximum degree six, it can be partitioned into seven independent sets and thus, since it has exactly $n-2$ vertices, it contains an independent set of size at least $(n-2)/7$ as required.
\end{proof}

In order to prove Theorem~\ref{T:D2} we need to improve upon the methods in the previous proof. Observe that the trivial bound $\Pr(A_y) \geqslant p^2$ holds with equality if and only if $y$ and $yx^{-1}$ are elements of order two in $G$. This might of course be the case for some $y$ but in this case, the number of neighbours of $y$ in the dependency graph $H$ will be less than six. (In fact it will be either one or two depending on whether $y$ commutes with $x$ or not.) So the idea is that either we can gain by improving the bound on $\Pr(A_y)$ for some $y$ or gain by showing that $y$ has fewer neighbours in the dependency graph. (And therefore the dependency graph has fewer edges and thus a larger independent set.) In the proof of Theorem~\ref{T:D2} we need to carefully control how much gain of each type we get for each particular vertex $y$.

\begin{proof}[Proof of Theorem~\ref{T:D2}]
As in the proof of Theorem~\ref{T:D2-worse constant} we aim to obtain an upper bound for the probability
\[ \Pr \left( \bigcap_{y \neq 1,x} \overline{A_y} \right), \]
where 1 denotes the identity element of $G$, $x$ is an element of $G$ different from the identity and $A_y$ is the event that the edges between $1$ and $y$ and between $x$ and $y$ both appear in $\Gamma$. 

As before, we let $H$ be the dependency graph on $G \setminus \{1,x\}$ for which two distinct vertices $y,z$ are adjacent if and only if 
\[
\{y,y^{-1},yx^{-1},xy^{-1}\} \cap \{z,z^{-1},zx^{-1},xz^{-1}\} \neq \emptyset.
\]
In particular, they are adjacent if and only if $z \in \{y^{-1},yx^{-1},xy^{-1},xy^{-1}x,yx,y^{-1}x\}$.

We divide the set of vertices of $G \setminus \{1,x\}$ into five types as follows:
\begin{itemize}
\item[Type~1:] All vertices $y$ which satisfy $y = y^{-1}$ and $xy^{-1}=yx^{-1}$. Observe that $A_y$ occurs if and only if $y$ and $xy^{-1}$ both appear in the generating set. In particular, since $xy^{-1} = xy \neq y$, we have that $\Pr(\overline{A_y}) = 1-p^2$. Note also that $N_H(y) = \{xy,yx\}$. (We are not excluding the possibility that $xy = yx$. We could have divide the class further into two subtypes according to whether $y$ commutes with $x$ or not but it turns out that this extra division is not needed.)
\item[Type~2a:] All vertices $y$ which satisfy $y = y^{-1}$ and $xy^{-1} \neq yx^{-1}$. Then $\Pr(\overline{A_y}) = 1-2p^2 + p^3$ and $N_H(y) = \{yx^{-1},xy,xyx,yx\}$.
\item[Type~2b:] All vertices $y$ which satisfy $y \neq y^{-1}$ and $xy^{-1} = yx^{-1}$. Then $\Pr(\overline{A_y}) = 1-2p^2 + p^3$ and $N_H(y) = \{y^{-1},xy^{-1},yx,y^{-1}x\}$.
\item[Type~3:] All vertices $y$ which satisfy $y = xy^{-1}$. Note that in this case we also have $y^{-1} = yx^{-1}$ and $x = y^2$. In particular we get $\Pr(\overline{A_y}) = 1 - 2p + p^2$ and $N_H(y) = \{y^{-1},y^3\}$. 
\item[Type~4:] All vertices $y$ which are not of the previous types. Note that in this case, all of $y,y^{-1},xy^{-1}$ and $yx^{-1}$ are distinct. In particular we get $\Pr(\overline{A_y}) = 1-4p^2 + 4p^3 - p^4$ and $N_H(y) = \{y^{-1},yx^{-1},xy^{-1},xy^{-1}x,yx,y^{-1}x\}$. 
\end{itemize}

We now prove two properties about $H$ that we will use.

\smallskip

\noindent
\textbf{Claim~1.} \emph{Vertices of Type~1 are adjacent only to vertices of Type~1 in $H$.}

\smallskip

To prove this we need to show that if $y = y^{-1}$ and $xy^{-1} = yx^{-1}$, then $z = z^{-1}$ and $xz^{-1} = zx^{-1}$ whenever $z = xy$ or $z = yx$. It is a simple check to see that this is indeed the case. This completes the proof of the claim.

\smallskip

\noindent
\textbf{Claim~2.} \emph{The subgraph of $H$ induced by vertices of Type~2a and Type~2b is four-colourable}

\smallskip

Since every vertex of Type~2a or Type~2b has degree at most four in $H$, by Brooks' theorem it is enough to check that this induced subgraph contains no clique on five vertices. Let us suppose that it does and let $y$ be one of its vertices. Suppose first that $y$ is of Type~2b. Then its neighbour  $xy^{-1}$ is of Type~2a. Indeed, to see this, we need to check that if $y \neq y^{-1}$ and $xy^{-1} = yx^{-1}$, then $xy^{-1} = (xy^{-1})^{-1}$ and $x(xy^{-1})^{-1} \neq (xy^{-1})x^{-1}$ which holds. So any clique on five vertices in this induced subgraph must contain a vertex, say $z$ of Type~2a. For this to happen, all of $xz,zx,zx^{-1}$ and $xzx$ must be distinct and either of Type~2a or of Type~2b. 

We claim now that $xzx$ is of Type~2a. To see this, observe that since $x(xzx)^{-1} = z^{-1}x^{-1} = zx^{-1}$ and $(xzx)x^{-1} = xz = xz^{-1}$, then $x(xzx)^{-1} \neq (xzx)x^{-1}$ and so $xzx$ cannot be of Type 2b. 

Therefore $(xzx) = (xzx)^{-1}$ and so $z = x^2zx^2$. Moreover the neighbours of $xzx$ in $H$ must be $x(xzx) = x^2zx = zx^{-1}, (xzx)x = x^{-1}z, (xzx)x^{-1} = xz$ and $x(xzx)x = z$. But we already know that the neighbours of $xzx$ in $H$ must be $z,xz,zx$ and $zx^{-1}$. It follows that $x^{-1}z = zx$, and thus $zx^{-1} = x^2zx = xz = xz^{-1}$, contradicting the fact that $z$ is of Type~2a. This completes the proof of the claim.

\smallskip

Now let $B_1$ be the set of all vertices of Type~1. Let $B_3$ be any maximal independent set containing vertices of Type~3, $B_4$ any maximal independent set containing vertices of Type~4 not adjacent to vertices in $B_3$ and finally $B_2$ any maximum independent set containing vertices of Types~2a and~2b not adjacent to any vertex of $B_3 \cup B_4$. Let us write $b_i = |B_i|$. Then by Claim~1 and the way $B_1,\ldots,B_4$ are define we have that
\begin{equation}\label{EQ:D2}
\Pr \left( \bigcap_{y \neq 1,x} \overline{A_y} \right) \leqslant  (1 - 2p^2 + p^3)^{b_2}(1-2p+p^2)^{b_3}(1-4p^2 + 4p^3 - p^4)^{b_4}\Pr \left( \bigcap_{y \in B_1} \overline{A_y} \right)
\end{equation}

\smallskip

\noindent
\textbf{Claim~3.} 
$b_1 + 4b_2 + 3b_3 + 7b_4 \geqslant n-2$
\smallskip

Since the vertices of Type 1 are adjacent only to vertices of Type 1, and since vertices of Type 3 have degree at most 2 and vertices of Type 4 have degree at most 6, then the number of vertices of Types 2a and 2b not incident to any vertex of $B_3 \cup B_4$ are at least $(n-2) - b_1 - 3b_3 - 7b_4$. But since by Claim~2 any induced subgraph of vertices of Types 2a and 2b is four-colourable, we must have $b_2 \geqslant \frac{1}{4}((n-2) - b_1 - 3b_3 - 7b_4)$ from which our claim follows. 

\smallskip

\noindent
\textbf{Claim~4.} 
$\Pr \left( \bigcap_{y \in B_1} \overline{A_y} \right) \leqslant (1 - p^2)^{b_1/2}$
\smallskip

Suppose that $x$ has order $k$. Then any $y$ which is of Type 1 is adjacent to $yx$ and to $xy = yx^{-1}$. It follows that any component of the dependency graph consisting of vertices of Type 1 has vertex set of the form $\{y,yx,yx^2,\ldots,yx^{k-1}\}$. Moreover, the event $A_{yx^{i}}$ happens if and only if both $yx^{i}$ and $yx^{i-1}$ appear in the generating set. Thus, the probability that none of the events $A_{yx^{i}}$ happens is the probability that when choosing elements from the set $\{1,2,\ldots,k\}$ independently at random with probability $p$, no two consecutive numbers are chosen, where $1$ and $k$ are considered consecutive.

If $k$ is even, then the claim follows immediately since the required probability is at most the probability that for each $1 \leqslant i \leqslant k/2$ not both of $2i-1$ and $2i$ appear, which is equal to $(1 - p^2)^{k/2}$.  

It remains to consider the case that $k \geqslant 3$ is odd. But then the required probability is at most the probability that for each $1 \leqslant i \leqslant (k-3)/2$ not both of $2i-1$ and $2i$ appear and moreover neither both of $k-2$ and $k-1$ appear and neither both of $k-1$ and $k$ appear. But this probability is equal to
\[
(1 - p^2)^{\frac{k-3}{2}}((1-p)^3 + 3p(1-p^2) + p^2(1-p)) = (1-p^2)^{{\frac{k-3}{2}}}(1 - 2p^2 + p^3) \leqslant (1 - p^2)^{k/2}.
\]
Here, we used the fact that $p^3 \leqslant p^2/2$ if $n$ is large enough .and also Bernoulli's inequality which states that $(1-x)^k \geqslant 1 - kx$ for $0 < x < 1$. This completes the proof of the claim.


\medskip

To complete the proof of the theorem, we use Equation~\eqref{EQ:D2}, Claims~3 and~4 and the fact that $(1 + x) \leqslant e^x$ for every real number $x$ to deduce that
\begin{align*}
\Pr \left( \bigcap_{y \neq 1,x} \overline{A_y} \right) &\leqslant \exp{\left\{-p^2\left(\frac{b_1}{2} + (2-p)b_2 + \left(\frac{2}{p} - 1\right)b_3 + (4 - 4p + p^2)b_4 \right)\right\}} \\
&\leqslant \exp{\left\{-\frac{p^2}{2}(b_1 + 4b_2 + 3b_3 + 7b_4) + p^3(b_2 + 4b_4) - pb_3(2 - 5p/2)\right\}} \\
&\leqslant \exp{\{-p^2(n-2)/2 + 4np^3\}} = o(1/n).
\end{align*}
Form the union bound it now follows that with high probability every vertex of $G$ is at distance at most two from vertex $1$ and since $G$ is vertex transitive it follows that with high probability $G$ has diameter at most 2.
\end{proof}

If we have a specific family of groups in mind,  it is sometimes possible to use the structure of the groups to modify the above proof and obtain better bounds. We proceed to do this for the family of cyclic groups. Even though this family is simple enough so that one can give a more direct proof we will proceed along the lines of the previous proof. In this proof we will make use of Janson's inequality as well, although it is not really needed. In this way, we prepare the ground for the proof of Theorem~\ref{E:D=2} which will follow in the next section.

\begin{proof}[Proof of Theorem~\ref{E:D>=3-Abelian}]
We denote the identity element of $G$ by 1, let $x \neq 1$ be an element of $G$ and for each $y \neq 1,x$ we write $A_y$ for the event that $\Gamma$ contains the edges $1y$ and $yx$. As before we aim to obtain an upper bound for the probability  
\[ \Pr \left( \bigcap_{y \neq 1,x} \overline{A_y} \right). \]

In fact, for reasons that will become clearer later, from now on we will further assume that $x^2 \neq 1$. Since $G$ is cyclic, there is at most one non-identity element $x$ with $x^2 = 1$ and for this we will use the fact that $\Pr(d(1,x) > 2) = o(1)$. This follows directly from the arguments in the proofs of Theorems~\ref{T:D2-worse constant} or~\ref{T:D2}.

Using the notation of the previous proof an element $y$ of $G$ would be of Type~4 (with respect to $x$) unless $y^2=1$, or $y^2 = x^2$, or $y^2 = x$ and there are at most six such possible choices for $y$ (including the choices $y=1$ and $y=x$). So letting $S$ denote the set of all $y$'s which do not satisfy any of the above equalities we have  
\[ \Pr \left( \bigcap_{y \neq 1,x} \overline{A_y} \right) \leqslant \Pr \left( \bigcap_{y \in S} \overline{A_y} \right)  \]
and it is enough to find an upper bound for the right hand side of the above equation. 

We now introduce some new notation in order to comply with the notation in our statement of Janson's inequality. We define an equivalence relation on $G$ by letting $y$ be equivalent to $z$ if and only if $y=z$ or $y=z^{-1}$. We write $[y]$ for the equivalence class of $y$ and we let $\Omega$ be the set of all equivalence classes. We let $R$ be a random subset of $\Omega$ where each equivalence class is chosen independently with probability $p$ or $2p-p^2$ depending on whether the equivalence class contains one or two elements. For each element $i$ of $G$ we let $F_i = \{[i],[xi^{-1}]\}$ and let $E_i$ be the event that $F_i \subseteq R$. Observe that for each $i$, $F_i = F_{xi^{-1}}$. So letting $I$ be a maximal subset of $S$ such that $i \in I \Rightarrow xi^{-1} \notin I$, then 
\[ 
\Pr \left( \bigcap_{y \in S} \overline{A_y} \right) = \Pr \left( \bigcap_{i \in I} \overline{E_i} \right).
\]
Observe that for each $i \in I$, since $i^2 \neq 1$ and $i^2 \neq x^2$ then $\Pr([i] \in R) = \Pr([xi^{-1}] \in R) = (2p-p^2)$. Since also $i^2 \neq x$, then $[i] \neq [xi^{-1}]$ and so $\Pr(E_i) = (2p-p^2)^2$. Finally, observe that for $i,j \in I$ with $i \neq j$, we have $F_i \cap F_j \neq \emptyset$ if and only if $j \in \{i^{-1},x^2i^{-1},xi,x^{-1}i\}$. (The case $j = xi^{-1}$ is excluded since by its definition, $I$ does not contain both $i$ and $xi^{-1}$.) We now claim that for every $i,j \in I$ with $i\neq j$ and $F_i \cap F_j \neq \emptyset$ it holds that $|F_i \cup F_j| = 3$. We proceed to verify this.

\begin{itemize}

\item[(a)] If $j = i^{-1}$ or $j = xi$ then $F_i \cup F_j =\{[i],[xi^{-1}],[xi]\}$. We claim that $|F_i \cup F_j| = 3$. Indeed, we already know that $[i] \neq [xi^{-1}]$. Since $x \neq 1$ and $i^2 \neq x$ we have $[i] \neq [xi]$. Moreover, since $i^2 \neq 1$ and $x^2 \neq 1$ we also have $[xi^{-1}] \neq [xi]$. 

\item[(b)] If $j = x^{-1}i$ or $j = x^{2}i^{-1}$, then $F_i \cup F_j =\{[i],[xi^{-1}],[x^2i^{-1}]\}$. We claim that $|F_i \cup F_j| = 3$. Indeed, we already know that $[i] \neq [xi^{-1}]$. If $j = x^{-1}i$, then $[xi^{-1}] = [j] \neq [xj^{-1}] = [x^2i{-1}]$, while if $j = x^{2}i^{-1}$, then $[x^2i^{-1}] = [j] \neq [xj^{-1}] = [xi^{-1}]$. Finally, since $i^2 \neq x^2$ and $x^2 \neq 1$ we also have $[i] \neq [x^2i^{-1}]$.
\end{itemize}

So for $i,j \in I$ with $i \neq j$ and $F_i \cap F_j \neq \emptyset$ we have that $\Pr(E_i \cap E_j) = (2p-p^2)^3$. Thus by Janson's inequality we get 
\begin{align*}
\Pr \left( \bigcap_{i \in I} \overline{E_i} \right) &\leqslant \exp\left\{-|I|(2p-p^2)^2 + 4|I|(2p-p^2)^3 \right\} \\
&= \exp\left\{-4p^2(1 + O(p))|I|\right\}.
\end{align*}
Finally, since $|I| \geqslant |S|/2 \geqslant (n-6)/2$ and $2p^2n \geqslant (1 + 2\eps)\log{n}$ we get that 
\[ \Pr(d(1,x) > 2) < \Pr \left( \bigcap_{i \in I} \overline{E_i} \right) = o(1/n).\]
Since this holds for every $x$ with at most one exception for which $\Pr(d(1,x) > 2) = o(1)$, the result follows from the union bound. 
\end{proof}

\section{Proof of Theorem~\ref{E:D=2}}

As in the proof of Theorem~\ref{E:D>=3-Abelian} we denote the identity element of $G$ by 1, we let $x \neq 1$ be an element of $G$ and for each $y \neq 1,x$ we write $A_y$ for the event that $\Gamma$ contains the edges $1y$ and $yx$. As before we aim to obtain an upper bound for the probability  
\[ \Pr \left( \bigcap_{y \neq 1,x} \overline{A_y} \right). \]
 
We begin by assuming that $x^2 \neq 1$ and $\cl(x) = \omega(1)$. The case in which either $x^2 = 1$ or $\cl(x) = O(1)$ is more complicated and will be treated later.

Using the notation of the previous proofs an element $y$ of $G$ would be of Type~4 (with respect to $x$) unless $y^2=1$, or $y^2 = x^2$, or $y^2 = x$. For reasons that will become clearer later, we let $S$ denote the set of all $y$'s such that $y^2,(y^{-1}x)^2 \notin \{1,x,x^{-1},x^2\}$. In particular, every element of $S$ is of Type~4. Furthermore, it follows from  Theorem~\ref{main-rep} that $S$ contains all but $o(n)$ elements of $G$. So 
\[ \Pr \left( \bigcap_{y \neq 1,x} \overline{A_y} \right) \leqslant \Pr \left( \bigcap_{y \in S} \overline{A_y} \right)  \]
and it is enough to find an upper bound for the right hand side of the above equation. 
We now recall the notation introduced in the proof of Theorem~\ref{E:D>=3-Abelian}. We define an equivalence relation on $G$ by letting $y$ be equivalent to $z$ if and only if $y=z$ or $y=z^{-1}$. We write $[y]$ for the equivalence class of $y$ and we let $\Omega$ be the set of all equivalence classes. We let $R$ be a random subset of $\Omega$ where each equivalence class is chosen independently with probability $p$ or $2p-p^2$ depending on whether the equivalence class contains one or two elements. For each element $i$ of $G$ we let $F_i = \{[i],[xi^{-1}]\}$ and let $E_i$ be the event that $F_i \subseteq R$. In contrast with the case of cyclic groups, we do not in general have $F_i = F_{xi^{-1}}$. Depending on the choice of $i,x$ we may have some other relations. More specifically,
\begin{itemize}
\item[(1)] If $x^2 = 1$ then $F_i = F_{ix}$.
\item[(2)] If $i^{-1}xi = x$ then $F_i = F_{xi^{-1}}$.
\item[(3)] If $i^{-1}xi = x^{-1}$ then $F_i = F_{i^{-1}}$.
\end{itemize}
For these reasons, we define $I=I(x)$ to be the largest subset of $S$ such that 
\begin{itemize}
\item[(1)] If $x^2 = 1$ then $i \in I \Rightarrow ix \notin I$.
\item[(2)] If $i^{-1}xi = x$ then $i \in I \Rightarrow xi^{-1} \notin I$.
\item[(3)] If $i^{-1}xi = x^{-1}$ then $i \in I \Rightarrow i^{-1} \notin I$.
\end{itemize}
It is important to note for later use that the above conditions are symmetric. For example, in condition (2), if we set $j = xi^{-1}$, then $j^{-1}xj = x$ and $i = xj^{-1}$. 

Observe that for each $i \in I$, since $i^2 \neq 1$ and $i^2 \neq x^2$ then $\Pr([i] \in R) = \Pr([xi^{-1}] \in R) = (2p-p^2)$. Since also $i^2 \neq x$, then $[i] \neq [xi^{-1}]$ and so $\Pr(E_i) = (2p-p^2)^2$. Finally, observe that for $i,j \in I$ with $i \neq j$, we have $F_i \cap F_j \neq \emptyset$ only if $j \in \{i^{-1},xi^{-1},ix^{-1},i^{-1}x,ix,xi^{-1}x\}$. We now claim that for every $i,j \in I$ with $i\neq j$ and $F_i \cap F_j \neq \emptyset$ it holds that $|F_i \cup F_j| = 3$. Indeed, we need to check the following cases:

\begin{itemize}

\item[(a)] If $j = i^{-1}$ then $F_i \cup F_j = \{[i],[xi^{-1}],[xi]\}$. We already know that $[i] \neq [xi^{-1}]$ and $[xi] = [xj^{-1}] \neq [j] = [i]$. So if $|F_i \cup F_j| \neq 3$ then $[xi^{-1}] = [xi]$ which can happen only if $xi^{-1} = xi$ or $ix^{-1} = xi$. None of these conditions hold as the first one implies that $i^2 = 1$ and so $i \notin I$, while the second one implies that $i^{-1}xi = x^{-1}$ and so by (3), we cannot have both $i,j \in I$. 

\item[(b)] If $j = xi^{-1}$ then $F_i \cup F_j = \{[i],[xi^{-1}],[xix^{-1}]\}$. We already know that $[i] \neq [xi^{-1}]$ and $[xi^{-1}] = [j] \neq [xj^{-1}] = [xix^{-1}]$. So if $|F_i \cup F_j| \neq 3$ then then $[i] = [xix^{-1}]$ which implies that $xi = ix$ or $xix^{-1} = i^{-1}$. None of these conditions hold as the first one implies that $i^{-1}xi = x$ and so by (2), we cannot have both $i,j \in I$, while the second one implies that $j^2 = (i^{-1}x)^2 = x^2$ and so $j \notin I$.
 
\item[(c)] If $j = ix^{-1}$ then $F_i \cup F_j = \{[i],[xi^{-1}],[x^2i^{-1}]\}$. We already know that $[i] \neq [xi^{-1}]$ and $[xi^{-1}] = [j] \neq [xj^{-1}] = [x^2i^{-1}]$. So if $|F_i \cup F_j| \neq 3$ then $[i] = [x^2i^{-1}]$ which implies that either $x^2 = i^2$ or $x^2 = 1$. None of these conditions hold as the first one implies that $i \notin I$, while the second one implies that $j = ix^{-1} = ix$ and so by (1), we cannot have both $i,j \in I$.

\item[(d)] If $j = i^{-1}x$ then $F_i \cup F_j = \{[i],[xi^{-1}],[i^{-1}]x\}$. We already know that $[i] \neq [xi^{-1}]$ and $[i^{-1}x] = [j] \neq [xj^{-1}] = [i]$. So if $|F_i \cup F_j| \neq 3$ then $[i^{-1}x] = [xi^{-1}]$ which implies that either $i^2 = 1$ or $xi = ix$. None of these conditions hold as the first one implies as the first one implies that $i \notin I$, while the second one implies that $i^{-1}xi = x$ and $j = xi^{-1}$ and so by (2), we cannot have both $i,j \in I$. 

\item[(e)] If $j = ix$ then $F_i \cup F_j = \{[i],[xi^{-1}],[ix]\}$. We already know that $[i] \neq [xi^{-1}]$ and $[ix] = [j] \neq [xj^{-1}] = [i]$. So if $|F_i \cup F_j| \neq 3$ then $[xi^{-1}] = [ix]$ which implies that either $xi^{-1} = ix$ or $x^2 = 1$. None of these conditions hold as the first one implies that $(i^{-1}x)^2 = x^2$ and so $i \notin I$, while the second condition by (1), implies that we cannot have both $i,j \in I$. 

\item[(f)] If $j = xi^{-1}x$ then $F_i \cup F_j = \{[i],[xi^{-1}],[xi^{-1}x]\}$. We already know that $[i] \neq [xi^{-1}]$ and $[xi^{-1}x] = [j] \neq [xj^{-1}] = [xi^{-1}]$. So if $|F_i \cup F_j| \neq 3$ then $[i] = [xi^{-1}x]$ which implies that either $i = xi^{-1}x$ or $i^{-1} = xi^{-1}x$. None of these conditions hold as the first one implies that $(i^{-1}x)^2 = 1$ and so $i \notin I$, while the second condition implies that $i^{-1}xi = x^{-1}$ and $j = i^{-1}$ and so by (3), we cannot have both $i,j \in I$. 
\end{itemize}

So applying Jensen's inequality, we obtain that 

\[ \Pr(d(1,x) > 2) \leqslant \Pr \left( \bigcap_{i \in I} \overline{E_i} \right) \leqslant \exp\left\{-4p^2(1 + O(p))|I|\right\}. \] 

To conclude, we need to estimate $|I(x)|$ from the nature of $x$. 

\smallskip

\noindent
\textbf{Case~I.} \emph{$x^2 \neq 1$ and $\cl(x) > 1/\eps$.}

\smallskip

By Theorem~\ref{commuting_pairs} there are at most $\eps n$ elements $i \in S$ such that $i^{-1}xi = x$. By the symmetry of condition (2) at most (in fact exaclty) half of these will not belong to $I(x)$. Similarly, there are at most $\eps n$ elements $i \in S$ such that $i^{-1}xi = x^{-1}$ and at most half of them do not belong to $I(x)$. Therefore $|I(x)| \geqslant (1 - \eps + o(1))n$ and so  $\Pr(d(1,x) > 2) \leqslant n^{-1 - 3\eps + 4\eps^2 + o(1)}$.

\smallskip

\noindent
\textbf{Case~II.} \emph{$x^2 = 1$ and $\cl(x) > 1/\eps$.}

\smallskip

As in Case~I, by removing at most $\eps n$ elements from $S$ we may assume that conditions (2) and (3) are satisfied. To satisfy condition (1), using its symmetry, we only need to remove at most half of the remaining elements. Therefore $|I(x)| \geqslant (1 - \eps + o(1))n/2$ and so $\Pr(d(1,x) > 2) \leqslant n^{-1/2 - 3\eps/2 + 2\eps^2 + o(1)}$.

\smallskip

\noindent
\textbf{Case~III.} \emph{$x^2 \neq 1$ and $\cl(x) \leqslant 1/\eps$.}

\smallskip

Observe that every $i\in S$ can satisfy at most one of the equalities $i^{-1}xi = x$ and $i^{-1}xi = x^{-1}$ (here we use the fact that $x^2 \neq 1$). So, by the symmetry of conditions (2) and (3) at most half of the elements of $S$ do not belong to $I$. Therefore $|I(x)| \geqslant (1 + o(1))n/2$ and so  $\Pr(d(1,x) > 2) \leqslant n^{-1/2 - 2\eps + o(1)}$.

\smallskip

\noindent
\textbf{Case~IV.} \emph{$x^2 = 1$ and $\cl(x) \leqslant 1/\eps$.}

\smallskip

As in the previous cases, the symmetry of conditions (1)-(3) guarantees that $|I| \geqslant |S|/4 = (1 + o(1))n/4$ and so  $\Pr(d(1,x) > 2) \leqslant n^{-1/4 - \eps + o(1)}$.

\smallskip

We can now use the union bound and the conditions (a)-(c) in the statement of the theorem to conclude that with high probability $d(1,x) > 2$ for every $x \neq 1$ and so with high probability $G$ has diameter at most 2.

\section{Proof of Theorem~\ref{Latin}}

Let us first assume that $p \geqslant \sqrt{\left(26 + \eps \right)\frac{\log{n}}{n}}$. We imitate the proof of Theorem~\ref{T:D2-worse constant}. There are two main differences in the proof. Firstly, because the graph might not be vertex transitive we need to prove that with high probability for every $x,y$ we have $d_{\Gamma}(x,y) \leqslant 2$. As in the proof of Theorem~\ref{T:D2-worse constant} we have that
\[ 
\Pr(d(x,y) > 2) \leqslant \Pr\left(\bigcap_{z \neq x,y}\overline{B_z}\right)
\]  
where $B_z = B_z(x,y)$ is the event that the edges between $x$ and $z$ and between $z$ and $y$ both appear in $\Gamma$. As before we have $\Pr(B_z) \geqslant p^2$. Moreover, the event depends only on the occurrence of the elements $L_{xz},L_{zx},L_{yz},L_{zy}$ in the generating set. So in the natural dependency graph $z$ will be adjacent to $w$ only if
\[
\{L_{xz},L_{zx},L_{yz},L_{zy}\} \cap \{L_{xw},L_{wx},L_{yw},L_{wy}\} \neq \emptyset
\]
and so by the definition of Latin squares it follows that the maximal dependency graph has maximum degree at most~12. A straightforward adaptation of the arguments in the proof of Theorem~\ref{T:D2-worse constant} show that 
\[
\Pr(d(x,y) > 2) \leqslant (1 - p^2)^{(n-2)/13}
\]
and so the result follows from the union bound.

For the other bound let us now assume that $p \leqslant \sqrt{\left(\frac{1}{4} + \eps \right)\frac{\log{n}}{n}}$. We imitate the proof of Theorem~\ref{T:D>=3}. Let us assume that the rows, columns and entries of $L$ are indexed by $\{1,2,\ldots,n\}$. It is enough to show that with high probability the number of $x \in \{2,\ldots,n\}$ which are at distance greater than two from vertex~1 is positive. As in the proof of Theorem~\ref{T:D>=3}, we define the graph $\Gamma_x$ on $\{1,2,\ldots,n\}$ where $i$ is adjacent to $j$ if and only if there is a $y$ such that $i \in \{L_{1y},L_{y1}\}$ and $j \in \{L_{xy},L_{yx}\}$ or vice versa. In particular, from the defining properties of a Latin square it follows that $\Gamma_x$ has maximum degree at most~8, minimum degree at least~1 and moreover, for every $i,j$ there are at most eight $\Gamma_x$'s in which $i$ and $j$ are adjacent. The rest of the proof is now identical to the proof of Theorem~\ref{T:D>=3} and is omitted.

\bigskip

\noindent
{\footnotesize
Demetres Christofides, Institute for Theoretical Computer Science, Faculty of Mathematics and Physics, \Malostranske\ \Namesti\ 25, 188 00 Prague, Czech Republic, \href{mailto:christofidesdemetres@gmail.com}{\tt christofidesdemetres@gmail.com}

\smallskip

\noindent Klas \Markstrom, Department of Mathematics and Mathematical Statistics, \Umea\ University, 90187 \Umea, Sweden, \href{mailto:klas.markstrom@math.umu.se}{\tt klas.markstrom@math.umu.se}
}

\end{document}